\documentclass[12 pt]{article}
\usepackage{amssymb, amsmath}
\usepackage{srcltx}

\begin{document}

\begin{center}
{\Large Isotropy of some quadratic forms and its applications on levels and
sublevels of algebras } \footnote{%
This paper was partially supported by the grant \textbf{UNESCO-UNITWIN
OCW/OER Initiative}, Handong Global University, South Korea.} 
\begin{equation*}
\end{equation*}%
Cristina Flaut

\bigskip

Faculty of Mathematics and Computer Science,

Ovidius University,

Bd. Mamaia 124, 900527, CONSTAN\c TA,

ROMANIA

e-mail: cflaut@univ-ovidius.ro; cristina\_flaut@yahoo.com

http://cristinaflaut.wikispaces.com/ 
\begin{equation*}
\end{equation*}
\end{center}

\textbf{Abstract.} {\small In this paper, we give some properties of the
levels and sublevels of algebras obtained by the Cayley-Dickson process. We
will emphasize how isotropy of some quadratic forms \ can influence the
levels and sublevels of algebras obtained by the Cayley-Dickson
process.\medskip }

{\small \ }

{\small Key Words: Cayley-Dickson process; Division algebra; Level and
sublevel of an algebra.}

{\small AMS Classification: 17A35, 17A20, 17A75, 17A45. } 
\begin{equation*}
\end{equation*}

\textbf{0.} \textbf{Introduction} 
\begin{equation*}
\end{equation*}

In [18], A. Pfister showed that if a field has a finite level this level is
a power of $\,\,2$ and any power of $2$ could be realised as the level of a
field. In the noncommutative case, the concept of level has many
generalisations. The level of division algebras is defined in the same
manner as for the fields.

In this paper, we give some properties of the levels and sublevels of
algebras obtained by the Cayley-Dickson process. We will emphasize how
isotropy of some quadratic forms \ can influence the levels and sublevels of
algebras obtained by the Cayley-Dickson process.\medskip \bigskip

\textbf{1. Preliminaries\bigskip }

In this paper, we assume that $K$ is a field and $charK\neq 2.$

For the basic terminology of quadratic and symmetric bilinear spaces, the
reader is referred to [22] or [11]. In this paper, we assume that all the
quadratic forms are nondegenerate.

A bilinear space $\left( V,b\right) $ \textit{represents} $\alpha \in K$ if
there is an element $x\in V,x\neq 0,$ with $b\left( x,x\right) =\alpha .$
The space is called \textit{universal} if $\,\left( V,b\right) $ represents
all $\alpha \in K.$ A vector $x\in V,x\neq 0$ is called \textit{isotropic}
if\thinspace \thinspace $b\left( x,x\right) =0,$ otherwise $x$ is called 
\textit{anisotropic}. If the bilinear space $\left( V,b\right) $, $V\neq
\{0\},$ contains an isotropic vector, then the space is called \textit{%
isotropic}. Every isotropic bilinear space $V$, $V\neq \{0\},$ is universal.
(See [22], Lemma 4.11., p.\thinspace 14])

The quadratic form $q:V\rightarrow K$ is called \textit{anisotropic} if $%
\,q\left( x\right) =0$ implies $x=0,$ otherwise $q$ is called \textit{%
isotropic}.

A quadratic form $\psi $ is a \textit{subform} of the form $\varphi $ if $%
\varphi \simeq \psi \perp \phi ,$ for some quadratic form $\phi .$ We denote 
$\psi <\varphi .$

Let $\varphi $ be a $n-$dimensional quadratic irreducible form over $K,$ $%
n\in N,n>1,$ which is not isometric to the hyperbolic plane$.$ We may
consider $\varphi $ as a homogeneous polynomial of degree $2$, 
\begin{equation*}
\varphi \left( X\right) =\varphi \left( X_{1},...X_{n}\right) =\sum
a_{ij}X_{i}X_{j},a_{ij}\in K^{\ast }.
\end{equation*}
The \textit{functions field of }$\varphi ,\,$denoted $K(\varphi ),$ is the
quotient field of the integral domain 
\begin{equation*}
K[X_{1},...,X_{n}]\,/\,\left( \varphi \left( X_{1},...,X_{n}\right) \right) .
\end{equation*}

Since $\left( X_{1},...,X_{n}\right) $ is a non-trivial zero, $\varphi $ is
isotropic over $K(\varphi ).\,$(See [22])\medskip

For $n\in \mathbb{N}-\{0\}\,$ a $n-$\textit{fold Pfister form} over $K$ is a
quadratic form of the type 
\begin{equation*}
<1,a_{1}>\otimes ...\otimes <1,a_{n}>,\text{ }a_{1},...,a_{n}\in K^{\ast }.
\end{equation*}%
A Pfister form is denoted by $\ll a_{1},a_{2},...,a_{n}\gg .$ For $n\in
N,n>1,\,\,\,$a Pfister form $\varphi $ can be written as 
\begin{equation*}
<1,a_{1}>\otimes ...\otimes
<1,a_{n}>=<1,a_{1},a_{2},...,a_{n},a_{1}a_{2},...,a_{1}a_{2}a_{3},...,a_{1}a_{2}...a_{n}>.
\end{equation*}%
If $\varphi =<1>\perp \varphi ^{\prime },$ then $\varphi ^{\prime }$ is
called \textit{the pure subform} of\thinspace $\varphi .$ A Pfister form is
hyperbolic if and only if is isotropic. This means that a Pfister form is
isotropic if and only if its pure subform is isotropic.( See [22] )

For the field $L$ it is defined%
\begin{equation*}
L^{\infty }=L\cup \{\infty \},
\end{equation*}
where $x+\infty =x$, for $x\in K,~x\infty =\infty $ for $x\in K^{\ast
},\infty \infty =\infty ,\frac{1}{\infty }=0,\frac{1}{0}=\infty $. 

An $L-$\textit{place }of the field $K$ is a map 
\begin{equation*}
\lambda :K\rightarrow L^{\infty }
\end{equation*}
with the properties:%
\begin{equation*}
\lambda \left( x+y\right) =\lambda \left( x\right) +\lambda \left( y\right)
,\lambda \left( xy\right) =\lambda \left( x\right) \lambda \left( y\right) ,
\end{equation*}%
whenever the right sides are defined.\bigskip 

\textbf{Theorem}([8], Theorem 3.3. )\textit{\ Let }$F$\textit{\ be a field
of characteristic }$\neq 2,$\textit{\ }$\varphi $\textit{\ be a quadratic
form over }$F$\textit{\ and }$\ K,L$\textit{\ extensions field of }$F.$%
\textit{\ If }$\varphi $\textit{\ }$_{K\text{ }}$\textit{\ is isotropic,
then there exist an }$F-$\textit{place from }$F\left( \varphi \right) $%
\textit{\ to }$K.\bigskip $

An algebra $A$ over $K$ is called \textit{quadratic} if $A$ is a unitary
algebra and, for all $x\in A,\,\,$there are $a,b\in K\,\,\,$ such that $%
x^{2}=ax+b1,$ $a,b\in K.$ The subset 
\begin{equation*}
A_{0}=\{x\in A-K\,\mid \,x^{2}\in K1\}
\end{equation*}
\thinspace \thinspace is a linear subspace of $A$ and 
\begin{equation*}
A=K\cdot 1\oplus A_{0}.
\end{equation*}

\textit{\ }A \textit{composition algebra} is an algebra $A$ \thinspace with
a non-degenerate quadratic form $q:A\rightarrow K,$ such that $q$ is
multiplicative, i.e. 
\begin{equation*}
q\left( xy\right) =q\left( x\right) q\left( y\right) ,\forall x,y\in A.
\end{equation*}
\thinspace A unitary composition algebra is called a \textit{Hurwitz algebra.%
} Hurwitz algebras have dimensions $1,2,4,8.$

Since over fields, the classical Cayley-Dickson process generates all
possible Hurwitz algebras, in the following, we briefly present the \textit{%
Cayley-Dickson process} and the properties of the algebras obtained.

Let $A$ be a finite dimensional unitary algebra over a field $\ K,$ with a 
\textit{scalar} \textit{involution} 
\begin{equation*}
\,\,\,\,\overline{\phantom{x}}:A\rightarrow A,a\rightarrow \overline{a},\,\,
\end{equation*}
i.e. a linear map satisfying the following relations:$\,\,\,$%
\begin{equation*}
\,\,\overline{ab}=\overline{b}\overline{a},\,\overline{\overline{a}}=a
\end{equation*}%
$\,\,$and 
\begin{equation*}
a+\overline{a},a\overline{a}\in K\cdot 1,
\end{equation*}
for all $a,b\in A$. The element $\,\overline{a}$ is called the \textit{%
conjugate} of the element $a.$ 

The linear form$\,\,$%
\begin{equation*}
\,\,t:A\rightarrow K\,,\,\,t\left( a\right) =a+\overline{a}
\end{equation*}%
is called the \textit{trace} of the element $a$ and \ the quadratic form 
\begin{equation*}
n:A\rightarrow K,\,\,n\left( a\right) =a\overline{a}
\end{equation*}%
is called the\textit{\ norm} of the element $a.$ It results that an algebra $%
A$ with a scalar involution is quadratic. $\,$If the quadratic form $n$ is 
\textit{anisotropic}, then the algebra $A$ is called \textit{anisotropic, }%
otherwise $A$ is \textit{isotropic.}

Let$\,\,\,\gamma \in K$ \thinspace be a fixed non-zero element. We define
the following algebra multiplication on the vector space $A\oplus A.$

\begin{equation}
\left( a_{1},a_{2}\right) \left( b_{1},b_{2}\right) =\left(
a_{1}b_{1}+\gamma \overline{b_{2}}a_{2},a_{2}\overline{b_{1}}%
+b_{2}a_{1}\right) .  \tag{1.1.}
\end{equation}%
We obtain an algebra structure over $A\oplus A.$ This algebra, denoted by $%
\left( A,\gamma \right) ,$ is called the \textit{algebra obtained from }$A$%
\textit{\ by the Cayley-Dickson process.} $A$ is canonically isomorphic with
the algebra 
\begin{equation*}
\,A^{\prime }=\{(a,0)\in A\oplus A\,\,\mid a\in A\},
\end{equation*}%
$\,$ where $\,A^{\prime }\,\,\,$ is a subalgebra of the algebra $\left(
A,\gamma \right) .\,\,\,$We denote $\left( 1,0\right) $ by $1$, where $(1,0)$
is the identity in $\left( A,\gamma \right) .\,$ Taking\thinspace \thinspace 
\begin{equation*}
\,\,u=\left( 0,1\right) \in A\oplus A,\,u^{2}=\gamma \cdot 1\in K\cdot 1,
\end{equation*}
it results that $\left( A,\gamma \right) =A\oplus Au$. We have $\dim \left(
A,\gamma \right) =2\dim A$.

Let $x\in \left( A,\gamma \right) ,x=\left( a_{1},a_{2}\right) .$ The map
\thinspace \thinspace \thinspace

\begin{equation*}
\overline{\phantom{x}}:\left( A,\gamma \right) \rightarrow \left( A,\gamma
\right) \,,\,\,\,\,x\rightarrow \bar{x}\,=\left( \overline{a}%
_{1},-a_{2}\right) ,
\end{equation*}%
is a scalar involution of the algebra $\left( A,\gamma \right) $, extending
the involution $\overline{\phantom{x}}\,\,\,$of the algebra $A$, therefore
the algebra $\left( A,\gamma \right) $ is quadratic. For $x\in \left(
A,\gamma \right) ,x=\left( a_{1},a_{2}\right) ,\,\,$we denote\thinspace
\thinspace 
\begin{equation*}
\,t\left( x\right) \cdot 1=x+\overline{x}=t(a_{1})\cdot 1\in K\cdot 1,
\end{equation*}%
\begin{equation*}
\,n\left( x\right) \cdot 1=x\overline{x}=(a_{1}\overline{a_{1}}-\gamma a_{2}%
\overline{a_{2}})\cdot 1=(n\left( a_{1}\right) -\gamma n(a_{2}))\cdot 1\in
K\cdot 1
\end{equation*}%
$\,$ and the scalars\thinspace 
\begin{equation*}
t\left( x\right) =t(a_{1}),\,\,\,n\left( x\right) =n\left( a_{1}\right)
-\gamma n(a_{2})\,
\end{equation*}%
$\,\,$are called the \textit{trace} and the \textit{norm} of the element $%
x\in $ $\left( A,\gamma \right) ,$ respectively.\thinspace \thinspace
It\thinspace \thinspace \thinspace follows\thinspace \thinspace \thinspace
that$\,$%
\begin{equation*}
\,x^{2}-t\left( x\right) x+n\left( x\right) =0,\,\forall x\in \left(
A,\gamma \right) .\,
\end{equation*}

\smallskip \thinspace If we take $A=K$ \thinspace and apply this process $t$
times, $t\geq 1,\,\,$we obtain an algebra over $K,\,\,A_{t}=K\{\alpha
_{1},...,\alpha _{t}\}.$ By induction, in this algebra we find a basis $%
\{1,f_{2},...,f_{q}\},q=2^{t},$ satisfying the properties:

\begin{equation*}
f_{i}^{2}=\alpha _{i}1,\,\,\alpha _{i}\in K,\alpha _{i}\neq 0,\,\,i=2,...,q.
\end{equation*}%
\begin{equation}
f_{i}f_{j}=-f_{j}f_{i}=\beta _{ij}f_{k},\,\,\beta _{ij}\in K,\,\,\beta
_{ij}\neq 0,i\neq j,i,j=\,\,2,...q,  \tag{1.2.}
\end{equation}%
$\beta _{ij}$ and $f_{k}$ being uniquely determined by $f_{i}$ and $f_{j}.$

If 
\begin{equation*}
\,x\in A_{t},x=x_{1}1+{\sum\limits_{i=2}^{q}}x_{i}f_{i},
\end{equation*}%
then 
\begin{equation*}
\bar{x}=x_{1}1-{\sum\limits_{i=2}^{q}}x_{i}f_{i}
\end{equation*}
and 
\begin{equation*}
t(x)=2x_{1},n\left( x\right) =x_{1}^{2}-{\sum\limits_{i=2}^{q}}\alpha
_{i}x_{i}^{2}.
\end{equation*}
$\,$In the above decomposition of $x$, we call $x_{1}$ the \textit{scalar
part} of $x$ and $x^{\prime \prime }=$ ${\sum\limits_{i=2}^{q}}x_{i}f_{i}$
the \textit{pure part} of $x.$ If we compute 
\begin{equation*}
x^{2}=x_{1}^{2}+x^{\prime \prime 2}+2x_{1}x^{\prime \prime }=
\end{equation*}%
\begin{equation*}
=x_{1}^{2}+\alpha _{1}x_{2}^{2}+\alpha _{2}x_{3}^{2}-\alpha _{1}\alpha
_{2}x_{4}^{2}+\alpha _{3}x_{5}^{2}-...-\left( -1\right) ^{t}({%
\prod\limits_{i=1}^{t}}\alpha _{i})x_{q}^{2}+2x_{1}x^{\prime \prime },
\end{equation*}
the scalar part of $x^{2}$ is represented by the quadratic form \newline
\begin{equation}
T_{C}=<1,\alpha _{1},\alpha _{2},-\alpha _{1}\alpha _{2},\alpha
_{3},...,\left( -1\right) ^{t}({\prod\limits_{i=1}^{t}}\alpha
_{i})>=<1,\beta _{2},...,\beta _{q}>  \tag{1.3.}
\end{equation}%
and, since 
\begin{equation*}
x^{\prime \prime 2}=\alpha _{1}x_{2}^{2}+\alpha _{2}x_{3}^{2}-\alpha
_{1}\alpha _{2}x_{4}^{2}+\alpha _{3}x_{5}^{2}-...-\left( -1\right) ^{t}({%
\prod\limits_{i=1}^{t}}\alpha _{i})x_{q}^{2}\in K,
\end{equation*}
it is represented by the quadratic form $T_{P}=T_{C}\mid
_{A_{0}}:A_{0}\rightarrow K,$%
\begin{equation}
T_{P}=<\alpha _{1},\alpha _{2},-\alpha _{1}\alpha _{2},\alpha
_{3},...,\left( -1\right) ^{t}({\prod\limits_{i=1}^{t}}\alpha _{i})>=<\beta
_{2},...,\beta _{q}>.  \tag{1.4.}
\end{equation}%
The quadratic form $T_{C}$ is called \textit{the trace form}, and $T_{P}$ 
\textit{the pure trace form} of the algebra $A_{t}.$ We remark that $%
T_{C}=<1>\perp T_{P},$ and the norm\newline
$n=n_{C}=<1>\perp -T_{P},$ resulting that\newline
\begin{equation*}
n_{C}=<1,-\alpha _{1},-\alpha _{2},\alpha _{1}\alpha _{2},\alpha
_{3},...,\left( -1\right) ^{t+1}({\prod\limits_{i=1}^{t}}\alpha
_{i})>=<1,-\beta _{2},...,-\beta _{q}>.\newline
\end{equation*}

The norm form $n_{C}$ has the form 
\begin{equation*}
n_{C}=<1,-\alpha _{1}>\otimes ...\otimes <1,-\alpha _{t}>
\end{equation*}
and it is a Pfister form.

Since the scalar part of any element $y\in A_{t}$ is $\frac{1}{2}t\left(
y\right) ,$ it follows that 
\begin{equation*}
T_{C}\left( x\right) =\frac{t\left( x^{2}\right) }{2}.\bigskip 
\end{equation*}

\textbf{Brown's construction of division algebras}\medskip \medskip

In 1967, R. B. Brown constructed, for every $t,$ a division algebra $A_{t}$
of dimension $2^{t}$ over the power-series field $K\{X_{1},X_{2},...,X_{t}%
\}. $ We will briefly demonstrate this construction, using polynomial rings
over $K$ and their field of fractions (the rational function field) instead
of power-series fields over $K$ (as it is done by R.B. Brown).

First of all, we remark that if an algebra $A$ is finite-dimensional, then
it is a division algebra if and only if $A$ does not contain zero divisors
(See [20]). For every $t$ we construct a division algebra $A_{t}$ over a
field $F_{t}.$ Let $X_{1},X_{2},...,X_{t}$ be $t$ algebraically independent
indeterminates over the field $K$ and $F_{t}=K\left(
X_{1},X_{2},...,X_{t}\right) $ be the rational function field$.$ For $%
i=1,...,t,$ we construct the algebra $A_{i}$ over the rational function
field $K\left( X_{1},X_{2},...,X_{i}\right) \,\,\,$by setting $\alpha
_{j}=X_{j}$ for $j=1,2,...,\,i.\,\,$Let $A_{0}=K.$ $\,$By\thinspace
\thinspace \thinspace \thinspace \thinspace induction over$\,\,\,i,$%
\thinspace \thinspace assuming that $A_{i-1}$ is a division algebra over the
field $F_{i-1}=K\left( X_{1},X_{2},...,X_{i-1}\right) $, we may
prove\thinspace \thinspace \thinspace \thinspace that\thinspace \thinspace
\thinspace \thinspace the\thinspace \thinspace \thinspace algebra$%
\,\,\,A_{i} $ is a division algebra over the field $F_{i}=K\left(
X_{1},X_{2},...,X_{i}\right) $.

$\,$Let 
\begin{equation*}
A_{F_{i}}^{i-1}=F_{i}\otimes _{F_{i-1}}A_{i-1}.
\end{equation*}
For $\alpha _{i}=X_{i}$ we apply the Cayley-Dickson process to algebra $%
A_{F_{i}}^{i-1}.$ The obtained algebra, denoted $A_{i},$ is an algebra over
the field $F_{i}$ and has dimension $2^{i}.$

Let $\,\,$%
\begin{equation*}
x=a+bv_{i},\,y=c+dv_{i},\,\,
\end{equation*}%
be nonzero elements in $A_{i}$ such that $xy=0,$ where $v_{i}^{2}=\alpha
_{i}.$ Since 
\begin{equation*}
xy=ac+X_{i}\bar{d}b+\left( b\bar{c}+da\right) v_{i}=0,
\end{equation*}%
we obtain 
\begin{equation}
ac+X_{i}\bar{d}b=0\,\,\text{\thinspace \thinspace }  \tag{2.1}
\end{equation}%
and 
\begin{equation}
b\bar{c}+da=0.  \tag{2.2.}
\end{equation}%
But, the elements $a,b,c,d$ $\in A_{F_{i}}^{i-1}\,\,$are non zero elements.
Indeed, we have:

i) If \thinspace $a=0$ and $b\neq 0,$ then $c=d=0\Rightarrow y=0,$ false;

ii)$\,$If $b=0$ and $a\neq 0,$ then $d=c=0\Rightarrow y=0,$ false;

iii)$\,$ If $c=0$ and\thinspace $\,d\neq 0,$ then $a=b=0\Rightarrow x=0,$
false;

iv) If\thinspace $d=0$\thinspace and $c\neq 0,$ then $a=b=0$ $\Rightarrow
x=0,\,$false.

This implies that $b\neq 0,a\neq 0,\,d\neq 0,\,c\neq 0.\,\,$If $%
\{1,f_{2},...,f_{2^{i-1}}\}$ is a basis in $A_{i-1},\,$then 
\begin{equation*}
a={{\sum\limits_{j=1}^{2^{i-1}}}}g_{j}(1\otimes f_{j})={{\sum%
\limits_{j=1}^{2^{i-1}}}}g_{j}f_{j},g_{j}\in F_{i},
\end{equation*}%
\begin{equation*}
g_{j}=\frac{g_{j}^{\prime }}{g_{j}^{\prime \prime }},g_{j}^{\prime
},g_{j}^{\prime \prime }\in K[X_{1},...,X_{i}],\,\,g_{j}^{\prime \prime
}\neq 0,\,j=1,2,...2^{i-1},
\end{equation*}
where \thinspace $K[X_{1},...,X_{t}]$ is the polynomial ring. Let $a_{2}$ $%
\,\,$be the less common multiple of\thinspace \thinspace \thinspace $%
g_{1}^{\prime \prime },....g_{2^{i-1}}^{\prime \prime },$ then we can write $%
a=\frac{a_{1}}{a_{2}},\,\,$where $a_{1}\in A_{F_{i}}^{i-1},a_{1}\neq 0.$
Analogously, $b=\frac{b_{1}}{b_{2}},c=\frac{c_{1}}{c_{2}},d=\frac{d_{1}}{%
d_{2}},b_{1},c_{1},d_{1}\in A_{F_{i}}^{i-1}-\{0\}$ and $%
a_{2},b_{2},c_{2},d_{2}\in K[X_{1},...,X_{t}]-\{0\}.$

If we replace in relations $\left( 2.1.\right) $ and $\left( 2.2.\right) ,$
we obtain 
\begin{equation}
a_{1}c_{1}d_{2}b_{2}+X_{i}\bar{d}_{1}b_{1}a_{2}c_{2}=0  \tag{2.3.}
\end{equation}%
and 
\begin{equation}
b_{1}\bar{c}_{1}d_{2}a_{2}+d_{1}a_{1}b_{2}c_{2}=0.  \tag{2.4.}
\end{equation}

If we denote $%
a_{3}=a_{1}b_{2},b_{3}=b_{1}a_{2},c_{3}=c_{1}d_{2},d_{3}=d_{1}c_{2},$ $%
a_{3},b_{3},c_{3},d_{3}\in A_{F_{i}}^{i-1}-\{0\},\,\,\,\,$relations $\left(
2.3.\right) $ and $\left( 2.4.\right) $ become

\begin{equation}
a_{3}c_{3}+X_{i}\bar{d}_{3}b_{3}=0  \tag{2.5.}
\end{equation}
and 
\begin{equation}
b_{3}\bar{c}_{3}+d_{3}a_{3}=0.  \tag{2.6.}
\end{equation}

Since the algebra $A_{F_{i}}^{i-1}=F_{i}\otimes _{F_{i-1}}A_{i-1}$ is an
algebra over $F_{i-1}$ with basis $X^{i}\otimes f_{j},$ $i\in \mathbb{N\,\,\,%
}$and $j=1,2,...2^{i-1},$ we can write $a_{3},b_{3},c_{3},d_{3}$ under the
form 
\begin{equation*}
a_{3}={\sum\limits_{j\geq m}}x_{j}X_{i}^{j},
\end{equation*}%
\begin{equation*}
b_{3}={\sum\limits_{j\geq n}}y_{j}X_{i}^{j},
\end{equation*}%
\begin{equation*}
c_{3}={\sum\limits_{j\geq p}}z_{j}X_{i}^{j},
\end{equation*}%
\begin{equation*}
d_{3}={\sum\limits_{j\geq r}}w_{j}X_{i}^{j},
\end{equation*}
where $x_{j},y_{j},z_{j},w_{j}\in A_{i-1},x_{m},y_{n},z_{p},w_{r}\neq 0.$
Since $A_{i-1}$ is a division algebra, we have $x_{m}z_{p}\neq
0,w_{r}y_{n}\neq 0,$ $y_{n}z_{p}\neq 0,w_{r}x_{m}\neq 0.$ Using relations $%
\left( 2.5.\right) $ and $\left( 2.6.\right) ,$ we have that $%
2m+p+r=2n+p+r+1,$ which is false. Therefore, the algebra $A_{i}$ is a
division algebra over the field $F_{i}=K\left( X_{1},X_{2},...,X_{i}\right)
\,\,\,$of dimension $2^{i}.$%
\begin{equation*}
\end{equation*}

\textbf{3.} \textbf{Main results}%
\begin{equation*}
\end{equation*}

The \textit{level} of the algebra $\,A,\,~$denoted by $s\left( A\right)
,\,\,\,$is the least integer $n$ such that $-1$ is a sum of $n$ squares in $%
A.$ The \textit{sublevel} of the algebra $A,\,\,$denoted by \underline{$s$}$%
\left( A\right) ,\,$ is the least integer $n$ such that $0$ is a sum of $%
\,\,n+1$ nonzero squares of elements in $A.$ If \thinspace these numbers do
not exist, then the level and sublevel are infinite. Obviously, \underline{$%
s $}$\left( A\right) \leq s\left( A\right) $.

Let $A$ be a division algebra over a field $K$ obtained by the
Cayley-Dickson process, of dimension $q=2^{t},T_{C},$ $T_{P}$ and $n_{C}$ be
its trace, pure trace and norm forms$.\medskip $

\textbf{Proposition 3.1.}\textit{\ With above notations, we have:}

\textit{i) If~~ }$s\left( A\right) \leq n$\textit{\ then }$-1$\textit{\ is
represented by the quadratic form }$n\times T_{C}.$

\textit{ii) }$-1$\textit{\ is a sum of }$~n$\textit{\ squares of pure
elements in }$A$\textit{\ if and only if the quadratic form }$n\times T_{P}$%
\textit{\ represents }$-1.$

\textit{iii)} \textit{For }$n\in \mathbb{N}-\{0\},$\textit{\ } \textit{if
the quadratic form }$<1>\perp n\times T_{P}\,\,\,$\textit{is isotropic over }%
$K,\,\,\,$\textit{then }$s(A)\leq n.\medskip $

\textbf{Proof.} i) Let $\;y\in
A,y=x_{1}+x_{2}f_{2}+...+x_{q}f_{q},\,x_{i}\in K,~$for all $i\in
\{1,2,...,q\}.$ \thinspace Using the notations given in the Introduction,$%
\,\,$we get 
\begin{equation*}
y^{2}=x_{1}^{2}+\beta _{2}x_{2}^{2}+...+\beta _{q}x_{q}^{2}+2x_{1}y^{\prime
\prime },
\end{equation*}
where 
\begin{equation*}
y^{\prime \prime }=x_{2}f_{2}+...+x_{q}f_{q}.
\end{equation*}
If $-1$ is a sum of $\,n\,\,$squares in $A,$ then%
\begin{equation*}
-1=y_{1}^{2}+...+y_{n}^{2}=
\end{equation*}%
\begin{equation*}
=\left( x_{11}^{2}\text{+}\beta _{2}x_{12}^{2}\text{+...+}\beta
_{q}x_{1q}^{2}\text{+}2x_{11}y_{1}^{\prime \prime }\right) \text{+...+}%
\left( x_{n1}^{2}\text{+}\beta _{2}x_{n2}^{2}\text{+...+}\beta _{q}x_{nq}^{2}%
\text{+}2x_{n1}y_{n}^{\prime \prime }\right) .
\end{equation*}
Then we have \newline
\begin{equation*}
-1={\sum\limits_{i=1}^{n}x_{i1}^{2}+\beta _{2}{\sum\limits_{i=1}^{n}}%
x_{i2}^{2}+...+\beta _{q}}{\sum\limits_{i=1}^{n}}x_{iq}^{2}
\end{equation*}
and\newline
\begin{equation*}
{\sum\limits_{i=1}^{n}}x_{i1}x_{i2}={\sum\limits_{i=1}^{n}}x_{i1}x_{i3}=...={%
\sum\limits_{i=1}^{n}}x_{i1}x_{in}=0,
\end{equation*}
then $n\times T_{C}$ represents $-1.$

ii) With the same notations, if $-1$ is a sum of $\,n\,\,$squares of pure
elements in $A,$ then 
\begin{equation*}
-1=y_{1}^{2}+...+y_{n}^{2}=
\end{equation*}%
\begin{equation*}
=\left( \beta _{2}x_{12}^{2}\text{+...+}\beta
_{q}x_{1q}^{2}+2x_{11}y_{1}^{\prime \prime }\right) \text{+...+}\left( \beta
_{2}x_{n2}^{2}\text{+...+}\beta _{q}x_{nq}^{2}+2x_{n1}y_{n}^{\prime \prime
}\right) .
\end{equation*}
We have $~$%
\begin{equation*}
-1={\beta _{2}{\sum\limits_{i=1}^{n}}x_{i2}^{2}+...+\beta _{q}}{%
\sum\limits_{i=1}^{n}}x_{iq}^{2}.
\end{equation*}
Therefore $n\times T_{P}$ ~represents~ $-1.$ Reciprocally, if \ $n\times
T_{P}$ ~represents~ $-1,$\ then 
\begin{equation*}
-1={\beta _{2}{\sum\limits_{i=1}^{n}}x_{i2}^{2}+...+\beta _{q}}{%
\sum\limits_{i=1}^{n}}x_{iq}^{2}.
\end{equation*}
Let 
\begin{equation*}
u_{i}={x_{i2}f}_{2}{+...+}x_{iq}^{2}f_{q}.
\end{equation*}%
$~$\ It results $t\left( u_{i}\right) =0$ and 
\begin{equation*}
u_{i}^{2}=-n\left( u_{i}\right) =\beta _{2}x_{i2}^{2}+...+\beta
_{q}x_{iq}^{2},
\end{equation*}
for all $i\in \{1,2,...,n\}.$ We obtain 
\begin{equation*}
-1=u_{1}^{2}+...+u_{n}^{2}.
\end{equation*}

iii) \textbf{Case 1.} If $-1\in K^{\ast 2},$ then $s\left( A\right) =1.$

\textbf{Case 2.} $-1\notin K^{\ast 2}.$ Since the quadratic form $<1>\perp
n\times T_{P}$ is isotropic then it is universal. It results that $<1>\perp
n\times T_{P}\,\,\,\,\,$represent $-1.$ Then, we have the elements $\alpha $ 
$\in K\,\,\,\,$ and \thinspace $p_{i}\in Skew(A),$ $i=1,...,n,\,\,\,$such
\thinspace that 
\begin{equation*}
-1=\alpha ^{2}+{\beta _{2}{\sum\limits_{i=1}^{n}}}p_{i2}^{2}+...+\beta _{q}{%
\sum\limits_{i=1}^{n}}p_{iq}^{2},
\end{equation*}
and not all of them are zero$.$

i) If\thinspace $\,\,\alpha =0,\,\,\,$then 
\begin{equation*}
-1={\beta _{2}{\sum\limits_{i=1}^{n}}}p_{i2}^{2}+...+\beta _{q}{%
\sum\limits_{i=1}^{n}}p_{iq}^{2}.
\end{equation*}
It \thinspace results \newline
\begin{equation*}
-1=\left( \beta _{2}p_{12}^{2}+...+\beta _{q}p_{1q}^{2}\right) +...+\left(
\beta _{2}p_{n2}^{2}+...+\beta _{q}p_{nq}^{2}\right) .
\end{equation*}
Denoting $\;$%
\begin{equation*}
u_{i}=p_{i2}f_{2}+...+p_{iq}f_{q},
\end{equation*}%
we have that $t\left( u_{i}\right) =0$ and 
\begin{equation*}
u_{i}^{2}=-n\left( u_{i}\right) =\beta _{2}p_{i2}^{2}+...+\beta
_{q}p_{iq}^{2},
\end{equation*}
for all $i\in \{1,2,...,n\}.$ We obtain $-1=u_{1}^{2}+...+u_{n}^{2}.$

ii) If $\alpha \neq 0,$ then $1+\alpha ^{2}\neq 0$ and 
\begin{equation*}
0=1+\alpha ^{2}+{\beta _{2}{\sum\limits_{i=1}^{n}}}p_{i2}^{2}+...+\beta _{q}{%
\sum\limits_{i=1}^{n}}p_{iq}^{2}.
\end{equation*}
Multiplying this relation with $1+\alpha ^{2}$ , it follows that 
\begin{equation*}
0=(1+\alpha ^{2})^{2}+{\beta _{2}{\sum\limits_{i=1}^{n}}}r_{i2}^{2}+...+%
\beta _{q}{\sum\limits_{i=1}^{n}}r_{iq}^{2}.
\end{equation*}
\thinspace\ Therefore%
\begin{equation*}
-1={\beta _{2}{\sum\limits_{i=1}^{n}}}r_{i2}^{\prime 2}+...+\beta _{q}{%
\sum\limits_{i=1}^{n}}r_{iq}^{\prime 2},
\end{equation*}
where 
\begin{equation*}
r_{ij}^{\prime }=r_{ij}(1+\alpha )^{-1},j\in \{2,3,...,q\}
\end{equation*}
and we apply case i). Therefore $s(A)\leq n.\Box \medskip $

\textbf{Proposition 3.2.} \textit{Let }$A$\textit{\ be a division algebra
obtained by the Cayley-Dickson process. The following statements are true:}

\textit{a)} \textit{If\thinspace \thinspace \thinspace }$n\in \mathbb{N}%
-\{0\},$\textit{\ } \textit{such that }$n=2^{k}-1,$ \textit{for} $k>1$, 
\textit{then }$s(A)\leq n$\textit{\ if and only if }$\,<1>\perp n\times
T_{P}\,$\textit{\ is isotropic.}

\textit{b) If }$\,-1$\textit{\ is a square in }$K,$\textit{\ then \thinspace 
}\underline{$s$}$\left( A\right) =s\left( A\right) =1.$

\textit{c) If\thinspace }$-1\notin K^{\ast 2},$\textit{\ then }$s\left(
A\right) =1$\textit{\ if and only if~ }$T_{C}$\textit{\ is isotropic.}%
\medskip

\textbf{Proof}. a) From Proposition 3.1, supposing that $s(A)\leq n,$ we have%
\begin{equation*}
-1={\sum\limits_{i=1}^{n}p_{i1}^{2}+\beta _{2}{\sum\limits_{i=1}^{n}}}%
p_{i2}^{2}+...+\beta _{q}{\sum\limits_{i=1}^{n}}p_{iq}^{2}
\end{equation*}
such that 
\begin{equation*}
\ {\sum\limits_{i=1}^{n}}p_{i1}p_{i2}={\sum\limits_{i=1}^{n}}%
p_{i1}p_{i3}=...={\sum\limits_{i=1}^{n}}p_{i1}p_{iq}=0.
\end{equation*}
For the level reasons, it results that 
\begin{equation*}
1+{\sum\limits_{i=1}^{n}p_{i1}^{2}\neq 0.}
\end{equation*}%
${\,}$\thinspace Putting\thinspace \thinspace \thinspace $p_{2^{k}1}=1$ and
\thinspace $p_{2^{k}2}=$ $p_{2^{k}3}=...$ $p_{2^{k}q}=0,$ we have 
\begin{equation}
0={\sum\limits_{i=1}^{n+1}p_{i1}^{2}+\beta _{2}{\sum\limits_{i=1}^{n+1}}}%
p_{i2}^{2}+...+\beta _{q}{\sum\limits_{i=1}^{n+1}}p_{iq}^{2}  \tag{3.1}
\end{equation}%
and%
\begin{equation*}
{\sum\limits_{i=1}^{n+1}}p_{i1}p_{i2}={\sum\limits_{i=1}^{n+1}}%
p_{i1}p_{i3}=...={\sum\limits_{i=1}^{n+1}}p_{i1}p_{iq}=0.
\end{equation*}
Multiplying $\left( 3.1.\right) $ by ${\sum\limits_{i=1}^{n+1}}p_{i1}^{2},$
since $\left( {\sum\limits_{i=1}^{n+1}}p_{i1}^{2}\right) ^{2}$is a square
and using Lemma from [22] p.151, for the products 
\begin{equation*}
{{\sum\limits_{i=1}^{n+1}}}p_{i2}^{2}{\sum\limits_{i=1}^{n+1}p_{i1}^{2},...,%
\sum\limits_{i=1}^{n+1}}p_{iq}^{2}{\sum\limits_{i=1}^{n+1}p_{i1}^{2}},
\end{equation*}%
we obtain 
\begin{equation}
0=\left( {\sum\limits_{i=1}^{n+1}}p_{i1}^{2}\right) ^{2}+\beta _{2}{%
\sum\limits_{i=1}^{n+1}}r_{i2}^{2}+...+\beta _{q}{\sum\limits_{i=1}^{n+1}}%
r_{iq}^{2},  \tag{3.2}
\end{equation}%
where 
\begin{equation*}
r_{i2},...r_{iq}\in K,n+1=2^{k}\,,
\end{equation*}%
$\,$ 
\begin{equation*}
r_{12}={\sum\limits_{i=1}^{n+1}}p_{i1}p_{i2}=0,\,r_{13}={\sum%
\limits_{i=1}^{n+1}}p_{i1}p_{i3}=0,...,r_{1q}={\sum\limits_{i=1}^{n+1}}%
p_{i1}p_{iq}=0.
\end{equation*}
Therefore, in the sums ${\sum\limits_{i=1}^{n+1}}r_{i2}^{2},\,...,{%
\sum\limits_{i=1}^{n+1}}r_{iq}^{2}\,\,\,$we have $n$ factors. From $\left(
3.2\right) ,$ we get that $<1>\perp n\times T_{P}$ is isotropic.

b) If $-1=a^{2}\in K\subset A,$ then \underline{$s$}$\left( A\right)
=s\left( A\right) =1.$

c) If $\,-1\notin K^{\ast 2}$ and $s\left( A\right) =1,$ then, there is an
element $y\in A\backslash K$ such that $-1=y^{2}.$ Hence $y\in A_{0},$ so $%
\overline{y}=-y.$ It results that 
\begin{equation*}
\left( 1+y\right) ^{2}=1+2y+y^{2}=2y
\end{equation*}
and 
\begin{equation*}
T_{C}\left( 1+y\right) =\frac{1}{2}t\left( \left( 1+y\right) ^{2}\right) =%
\frac{1}{2}\left( 2y+\overline{2y}\right) =y-y=0.
\end{equation*}
Therefore $T_{C}$ is isotropic.

Conversely, if $T_{C}$ is isotropic, then there is $y\in A,\,y\neq 0,$ such
that 
\begin{equation*}
T_{C}\left( y\right) =0=y_{1}^{2}+\beta _{2}y_{2}^{2}+...+\beta
_{q}y_{q}^{2}.
\end{equation*}
If $y_{1}=0,$ then $T_{C}\left( y\right) =T_{P}\left( y\right) =0,$ so $y=0,$
which is false.\thinspace \thinspace \thinspace If $y_{1}\neq 0,$ then 
\begin{equation*}
-1=\left( \left( \frac{y_{2}}{y_{1}}\right) f_{2}+...+\left( \frac{y_{q}}{%
y_{1}}\right) f_{q}\right) ^{2},
\end{equation*}
obtaining $s\left( A\right) =1.\Box \medskip $

\textbf{Proposition 3.3.} \textit{The quadratic form }$2^{k}\times T_{C}$%
\textit{\ is isotropic if and only if }$<1>\perp 2^{k}\times T_{P}$\textit{\
is isotropic.\medskip }

\textbf{Proof.} \thinspace Since the form $<1>\perp 2^{k}\times T_{P}$ is a
subform of the form $2^{k}\times T_{C}$, if the form $<1>\perp 2^{k}\times
T_{P}$ is isotropic, we have that \thinspace $2^{k}\times T_{C}$ \thinspace
is isotropic.

For the converse, supposing that $2^{k}\times T_{C}$ is isotropic, then we
get 
\begin{equation}
{\sum\limits_{i=1}^{2^{k}}}p_{i}^{2}+\beta _{2}{\sum\limits_{i=1}^{2^{k}}}%
p_{i2}^{2}+...+\beta _{q}{\sum\limits_{i=1}^{2^{k}}}p_{iq}^{2}=0,  \tag{3.3}
\end{equation}
where $p_{i},p_{ij}\in K,i=1,...,2^{k},\,j\in 2,...,q$ and \thinspace some
of the elements $p_{i}$ and\thinspace \thinspace $p_{ij}$ are nonzero$.$

If\thinspace \thinspace \thinspace \thinspace $p_{i}=0$, $\forall
i=1,...,2^{k},$ then $2^{k}\times T_{P}$ is isotropic, therefore $<1>\perp
2^{k}\times T_{P}$ is isotropic.

If 
\begin{equation*}
{\sum\limits_{i=1}^{2^{k}}}p_{i}^{2}\neq 0,
\end{equation*}%
then, multiplying relation $\left( 3.3\right) $ with ${\sum%
\limits_{i=1}^{2^{k}}}p_{i}^{2}$ and using Lemma from [22] p.151, for the
products 
\begin{equation*}
{\sum\limits_{i=1}^{2^{k}}}p_{i2}^{2}{\sum\limits_{i=1}^{2^{k}}}%
p_{i}^{2},...,{\sum\limits_{i=1}^{2^{k}}}p_{iq}^{2}{\sum\limits_{i=1}^{2^{k}}%
}p_{i}^{2},\,
\end{equation*}%
$\,$we obtain 
\begin{equation*}
\ ({\sum\limits_{i=1}^{2^{k}}}p_{i}^{2})^{2}+\beta _{2}{\sum%
\limits_{i=1}^{2^{k}}}r_{i2}^{2}+...+\beta _{q}{\sum\limits_{i=1}^{2^{k}}}%
r_{iq}^{2}=0,
\end{equation*}
then $<1>\perp 2^{k}\times T_{P}$ is isotropic.

For the level reason\medskip ,\thinspace \thinspace \thinspace the relation $%
{\sum\limits_{i=1}^{2^{k}}}p_{i}^{2}=0,$ for some $p_{i}\neq 0,$ does not
work. Indeed, supposing that $p_{1}\neq 0,$ we obtain 
\begin{equation*}
-1={\sum\limits_{i=2}^{2^{k}}(}p_{i}p_{1}^{-1})^{2},
\end{equation*}
false.$\Box $ $\medskip $

\textbf{Proposition 3.4.\thinspace \thinspace \thinspace \thinspace }\textit{%
\thinspace Let }$A$\textit{\ be an algebra over a field }$K$\textit{\
obtained by the Cayley-Dickson process, of dimension }$q=2^{t},T_{C}$\textit{%
\ and }$T_{P}$\textit{\ be its trace and pure trace forms}$.$\textit{\ If }$%
t\geq 2$ \thinspace \textit{and \ }$2^{k}\times T_{P}$\textit{\ is isotropic
over }$K$\textit{,\thinspace \thinspace }$k\geqslant 0,$\textit{\ then }$%
\left( 1+[\frac{2}{3}2^{k}]\right) \times T_{P}$\textit{\ is isotropic over }%
$K.\medskip $

\textbf{Proof. } If\thinspace \thinspace \thinspace $2^{k}\times T_{P}$ is
isotropic then $2^{k}\times -T_{P}$ is isotropic. Since 
\begin{equation*}
2^{k}\times n_{C}=2^{k}\times (<1>\bot -T_{P})
\end{equation*}
and $n_{C}$ is a Pfister form, we have $2^{k}\times n_{C}$ \thinspace is a
Pfister form. Since $2^{k}\times -T_{P}$ is a subform of $2^{k}\times n_{C},$
it results that $2^{k}\times n_{C}$ is isotropic, then it is hyperbolic.
Therefore 
\begin{equation*}
2^{k}\times n_{C}\simeq <1,1,...,1,-1,...,-1>
\end{equation*}%
(there are $2^{k+t-1}$ of $-1$ and $2^{k+t-1}$ of $1$)$.$ Multiplying by $-1,
$ we have that $2^{k}\times (<-1>\bot T_{P})$\thinspace \thinspace
\thinspace \thinspace \thinspace is hyperbolic, then has a totally isotropic
subspace of dimension $2^{k+t-1}$. It results that each subform of the form $%
2^{k}\times (<-1>\bot T_{P})$ of dimension greater or equal to $2^{k+t-1}$
is isotropic. Since 
\begin{equation*}
(2^{t}-1)(1+[\frac{2}{3}2^{k}])>(2^{t}-1)(\frac{2}{3}%
2^{k})>2^{t-1}2^{k}=2^{k+t-1},
\end{equation*}%
then $\left( 1+[\frac{2}{3}2^{k}]\right) \times T_{P}$\textit{\ }is
isotropic over\textit{\ }$K.\Box \medskip $

\textbf{Proposition 3.5.} \textit{Let }$A$\textit{\ be an algebra over a
field }$K$\textit{\ obtained by the Cayley-Dickson process, of dimension }$%
q=2^{t},T_{C}$\textit{\ and }$T_{P}$\textit{\ be its trace and pure trace
forms}$.\,\,$\textit{Let }$n=2^{k}-1.\mathit{\ \,}$\textit{If\thinspace
\thinspace \thinspace \thinspace }$t\geq 2\,\,$ \textit{and}\thinspace
\thinspace $k>1$ \textit{then} \underline{$s$}$\left( A\right) \leq 2^{k}-1$ 
\textit{if and only if} $\,\,<1>\bot (2^{k}-1)\times T_{P}$ \textit{is
isotropic}.\medskip

\textbf{Proof.} \ First, we prove the following result:

\textbf{Lemma. }\textit{For }$n=2^{k}-1,$ \underline{$s$}$\left( A\right)
\leq n$\textit{\ if and only if~\thinspace }$<1>\perp \left( n\times
T_{P}\right) $\textit{\ is isotropic or }$\left( n+1\right) \times T_{P}$%
\textit{\ is isotropic.\medskip }

\textbf{Proof of the Lemma.} Since \thinspace \underline{$s$}$\left(
A\right) \leq s\left( A\right) ,$ if\thinspace \thinspace $<1>\perp \left(
n\times T_{P}\right) $ is isotropic, then, from Proposition 3.2, we have 
\underline{$s$}$\left( A\right) \leq n.$ If $\left( n+1\right) \times T_{P}$
is isotropic, then there are the elements $p_{ij}\in K,i=1,...,2^{k},j\in
2,...,q,$ some of them are nonzero, such that {\thinspace \thinspace
\thinspace }%
\begin{equation*}
\beta _{2}{\sum\limits_{i=1}^{2^{k}}}p_{i2}^{2}+...+\beta _{q}{%
\sum\limits_{i=1}^{2^{k}}}p_{iq}^{2}=0.
\end{equation*}
We obtain 
\begin{equation*}
0=\left( \beta _{2}p_{12}^{2}+...+\beta _{q}p_{1q}^{2}\right) +...+\left(
\beta _{2}p_{n2}^{2}+...+\beta _{q}p_{nq}^{2}\right) .
\end{equation*}
Denoting $\;$%
\begin{equation*}
u_{i}=p_{i2}f_{2}+...+p_{iq}f_{q},
\end{equation*}
we have $t\left( u_{i}\right) =0$ and 
\begin{equation*}
u_{i}^{2}=-n\left( u_{i}^{2}\right) =\beta _{2}p_{i2}^{2}+...+\beta
_{q}p_{iq}^{2},
\end{equation*}
for all $i\in \{1,2,...,n\}.$ Therefore 
\begin{equation*}
0=u_{1}^{2}+...+u_{n}^{2}.
\end{equation*}
It results that \underline{$s$}$\left( A\right) \leq n.$

Conversely, if \underline{$s$}$\left( A\right) \leq n,$ then there are the
elements $y_{1},...,y_{n+1}\in A,$ some of them must be nonzero, such that 
\begin{equation*}
\,0=y_{1}^{2}+...+y_{n+1}^{2}.
\end{equation*}%
$\,\,$As in the proof of Proposition 3.1., we obtain 
\begin{equation*}
0={\sum\limits_{i=1}^{n+1}\mathit{x}_{i1}^{2}+\beta _{2}{\sum%
\limits_{i=1}^{n+1}}\mathit{x}_{i2}^{2}+...+\beta _{q}}{\sum%
\limits_{i=1}^{n+1}}x_{iq}^{2}
\end{equation*}
and 
\begin{equation*}
{\sum\limits_{i=1}^{n+1}}x_{i1}x_{i2}={\sum\limits_{i=1}^{n+1}}%
x_{i1}x_{i3}...={\sum\limits_{i=1}^{n+1}}x_{i1}x_{iq}=0.
\end{equation*}%
$\,\,$If all$\,\,\,x_{i1}=0,$ then $\left( n+1\right) \times T_{P}$ is
isotropic. If ${\sum\limits_{i=1}^{n+1}\mathit{x}_{i1}^{2}\neq 0,}$ then $%
\left( n+1\right) \times T_{C}$ is isotropic, or multiplying the last
relation with ${\sum\limits_{i=1}^{2^{k}}}x_{i1}^{2}$ \thinspace and using
Lemma from [22] p.151 for the products 
\begin{equation*}
{\sum\limits_{i=1}^{2^{k}}}x_{i2}^{2}{\sum\limits_{i=1}^{2^{k}}}%
x_{i1}^{2},...,{\sum\limits_{i=1}^{2^{k}}}x_{iq}^{2}{\sum%
\limits_{i=1}^{2^{k}}}x_{i1}^{2},
\end{equation*}
we obtain that $<1>\perp \left( n\times T_{P}\right) $ is isotropic. For
level reason of the field$,$ the relation ${\sum\limits_{i=1}^{n+1}\mathit{x}%
_{i1}^{2}=0}$ for some $x_{i1}\neq 0$ is false.

Using the above Lemma, we have that \underline{$s$}$\left( A\right) \leq
2^{k}-1$ if and only if $\,\,<1>\perp \left( n\times T_{P}\right) $\textit{\ 
}is isotropic or\textit{\ }$\left( n+1\right) \times T_{P}$\textit{\ }is
isotropic. In this case, we prove that $2^{k}\times T_{P}$ is isotropic
implies $<1>\bot (2^{k}-1)\times T_{P}$ is isotropic. If\thinspace
\thinspace $\,2^{k}\times T_{P}$ \thinspace \thinspace isotropic over $K$
then \thinspace $\left( 1+[\frac{2}{3}2^{k}]\right) \times T_{P}\,\,$ is
isotropic over $K.$ If\ \ $k\geq 2,$ then $\left( 1+[\frac{2}{3}%
2^{k}]\right) \leq 2^{k}-1$ and we have that $\left( 1+[\frac{2}{3}%
2^{k}]\right) \times T_{P}$ is an isotropic subform of the form $<1>\bot
(2^{k}-1)\times T_{P}$.$\Box \medskip $

\textbf{Remark 3.6.} Using the above notations, if the algebra $A$ is an
algebra obtained by the Cayley-Dickson process, of dimension greater than $2$
and if $\,n_{C\text{ }}$ is isotropic, then $s\left( A\right) =\underline{s}%
\left( A\right) =1.$ Indeed, if $-1$\textit{\ }is a square in\textit{\ }$K,$
the statement results from above. If $-1\notin K^{\ast 2},$ since $%
n_{C}=<1>\perp -T_{P}$ and $n_{C\text{ }}$ is a Pfister form, we obtain that 
$-T_{P}$ is isotropic, therefore $T_{C}$ is isotropic and, from above
proposition, we have that $s\left( A\right) =\underline{s}\left( A\right)
=1.\medskip $\medskip

\textbf{Proposition 3.7.} \textit{Let }$A$\textit{\ be an algebra over a
field }$K$\textit{\ obtained by the Cayley-Dickson process, of dimension }$%
q=2^{t},T_{C}$\textit{\ and }$T_{P}$\textit{\ be its trace and pure trace
forms}$.$ \textit{If \ }$k\geq t,$\textit{then }$s\left( A\right) \leq 2^{k}$%
\textit{\ if and only if the form }$\left( 2^{k}+1\right) \times <1>\perp
\left( 2^{k}-1\right) \times T_{P}$\textit{\ is isotropic}$.\medskip $

\textbf{Proof. }\ First, we prove the following result:\medskip

\textbf{Lemma.~} $s\left( A\right) \leq 2^{k}$\textit{\ }$=n$\textit{\ if
and only if~~}$\left( 2^{k}+1\right) \times <1>\perp \left( 2^{k}-1\right)
\times T_{P}~$\textit{\ or\ }$\ <1>\perp 2^{k}\times T_{P}$\textit{\ is
isotropic}$.\medskip $

\textbf{Proof of the Lemma.} If~ $s\left( A\right) \leq 2^{k},$ therefore $%
-1=y_{1}^{2}+...+y_{n}^{2},$ where $y_{i}\in A.~$\ Using above notations, we
obtain 
\begin{equation}
0=1+{\sum\limits_{i=1}^{n}x_{i1}^{2}+\beta _{2}{\sum\limits_{i=1}^{n}}%
x_{i2}^{2}+...+\beta _{q}}{\sum\limits_{i=1}^{n}}x_{iq}^{2}  \tag{*}
\end{equation}%
and 
\begin{equation*}
\ {\sum\limits_{i=1}^{n}}x_{i1}x_{i2}={\sum\limits_{i=1}^{n}}%
x_{i1}x_{i3}=...={\sum\limits_{i=1}^{n}}x_{i1}x_{iq}=0.
\end{equation*}%
If\ ${\sum\limits_{i=1}^{n}x_{i1}^{2}=0,}$ then $<1>\perp 2^{k}\times T_{P}.$
If~ $1+{\sum\limits_{i=1}^{n}x_{i1}^{2}=0,}$we have \ $2^{k}\times T_{P}$
isotropic, therefore $<1>\perp 2^{k}\times T_{P}$ isotropic. If \ $a={%
\sum\limits_{i=1}^{n}x_{i1}^{2}\neq 0,}$ multiplying relation $\left( \ast
\right) $ with \ ${\sum\limits_{i=1}^{n}x_{i1}^{2}}$ and using Lemma from
[22] p.151, for the products%
\begin{equation*}
{{\sum\limits_{i=1}^{n}}}x_{i2}^{2}{\sum\limits_{i=1}^{n}x_{i1}^{2},...,\sum%
\limits_{i=1}^{n}}x_{iq}^{2}{\sum\limits_{i=1}^{n}x_{i1}^{2},}\,
\end{equation*}%
$\,$we obtain 
\begin{equation*}
\ {\sum\limits_{i=1}^{n}x_{i1}^{2}+}a^{2}+\beta _{2}{\sum\limits_{i=1}^{n}}%
r_{i2}^{2}+...+\beta _{q}{\sum\limits_{i=1}^{n}}r_{iq}^{2}=0,
\end{equation*}
where 
\begin{equation*}
r_{i2},...r_{iq}\in K,
\end{equation*}%
$\,$ 
\begin{equation*}
r_{12}={\sum\limits_{i=1}^{n}}x_{i1}x_{i2}=0,\,r_{13}={\sum\limits_{i=1}^{n}}%
x_{i1}x_{i3}=0,...,r_{1q}={\sum\limits_{i=1}^{n}}p_{i1}p_{iq}=0.
\end{equation*}
Therefore, in the sums ${\sum\limits_{i=1}^{n-1}}r_{i2}^{2},\,...,{%
\sum\limits_{i=1}^{n-1}}r_{iq}^{2}\,\,\,$we have $n-1$ factors. It results 
\begin{equation*}
\left( 2^{k}+1\right) \times <1>\perp \left( 2^{k}-1\right) \times T_{P}
\end{equation*}
isotropic.

Conversely, if \ $<1>\perp 2^{k}\times T_{P}$\textit{\ }$\ $\ is isotropic,
from Proposition 3.1. iii), we have $s\left( A\right) \leq 2^{k}.$ If \ 
\textit{~}$\left( 2^{k}+1\right) \times <1>\perp \left( 2^{k}-1\right)
\times T_{P}$ is isotropic, then there is a nonzero element $\ a\in
D_{K}\left( 2^{k}\times <1>\right) $ such that $-a\in D_{K}\left( <1>\perp
\left( 2^{k}-1\right) \times T_{P}\right) .$ It results that 
\begin{equation*}
-a=b^{2}+{\beta _{2}B}_{2}{+...+\beta _{q}}B_{q},
\end{equation*}
where 
\begin{equation*}
b\in K,B_{2},...,B_{q}\in D_{K}\left( \left( 2^{k}-1\right) \times
<1>\right) \cup \{0\},
\end{equation*}
so that 
\begin{equation*}
-1=\frac{1}{a^{2}}(b^{2}a+{\beta _{2}B}_{2}a{+...+\beta _{q}}B_{q}a).
\end{equation*}
Supposing that $a={\sum\limits_{i=1}^{n}x_{i}^{2},}$ there are the elements $%
y_{i2}\in K$ such that ${\sum\limits_{i=1}^{n}y_{i2}^{2}=aB}_{2}$ and ${%
\sum\limits_{i=1}^{n}x_{i}y_{i2}=0.}$ Indeed, if $B_{2}=0,$ we put $%
y_{i2}=0,i\in \{1,2,...,n\}.$ If \ $B_{2}\neq 0,$ we have $<a,{aB}%
_{2}>\simeq a<1,{B}_{2}>.$ Since 
\begin{equation*}
B_{2}\in D_{K}\left( \left( 2^{k}-1\right) \times <1>\right) 
\end{equation*}
and 
\begin{equation*}
a\cdot 2^{k}\times <1>\simeq 2^{k}\times <1>,
\end{equation*}
it results that $a<1,{B}_{2}>$ is a subform of $2^{k}\times <1>,$ then
elements $y_{i2}$ exist. In the same way, we get the elements $y_{ij}\in
K,j\in \{3,...,q\}$ such that $aB_{j}={\sum\limits_{i=1}^{n}y_{ij}^{2}}$ and 
${\sum\limits_{i=1}^{n}x_{i}y_{ij}=0,}$ obtaining that 
\begin{equation*}
\sum\limits_{i=1}^{n}\frac{1}{a^{2}}\left(
ax_{i}+y_{i2}f_{2}+...+y_{iq}f_{q}\right) ^{2}=
\end{equation*}%
\begin{equation*}
=\frac{1}{a^{2}}(b^{2}a+{\beta _{2}B}_{2}a{+...+\beta _{q}}B_{q}a)=-1,
\end{equation*}
then $s\left( A\right) \leq n.\Box $

Using above Lemma, we have $\,\,$\textit{\ \ }$s\left( A\right) \leq 2^{k}$
if and only if $\left( 2^{k}+1\right) \times <1>\perp \left( 2^{k}-1\right)
\times T_{P}~$\textit{\ }or\textit{\ }$\ <1>\perp 2^{k}\times T_{P}$ is
isotropic. If $\ k\geq t,$ we prove that $<1>\perp 2^{k}\times T_{P}$
isotropic implies $\left( 2^{k}+1\right) \times <1>\perp \left(
2^{k}-1\right) \times T_{P}$ isotropic. Indeed, $<1>\perp 2^{k}\times T_{P}$
is isotropic if and only if $2^{k}\times T_{C}$ is isotropic, from
Proposition 3.3.. \ From [3] \ Proposition 1.4., it results that the Witt
index(i.e. the dimension of maximal totally isotropic subspace) of $%
2^{k}\times T_{C}$ is greater or equal with $2^{k}$. Therefore $2^{k}\times
T_{C}$ has a totally isotropic subspace of dimension $\geq 2^{k}.$ The form $%
2^{k}\times <1>\perp \left( 2^{k}-1\right) \times T_{P}$ is a subform of the
forms $2^{k}\times T_{C}$ and $\left( 2^{k}+1\right) \times <1>\perp \left(
2^{k}-1\right) \times T_{P}$ and has dimension $\left( 2^{t}-1\right) \left(
2^{k}-1\right) +2^{k}.$ Since $\left( 2^{t}-1\right) \left( 2^{k}-1\right)
+2^{k}=2^{k+t}-2^{t}+1>2^{k},$ for $k\geq t,$ we have that $2^{k}\times
<1>\perp \left( 2^{k}-1\right) \times T_{P}$ is isotropic, then $\left(
2^{k}+1\right) \times <1>\perp \left( 2^{k}-1\right) \times T_{P}$ is
isotropic.$\Box \medskip $

\textbf{Proposition 3.8.} \textit{Let }$K$ \textit{be a field.\ }

\textit{i) If\thinspace } $k\geq 2,$ \textit{then\ } \underline{$s$}$\left(
A\right) \leq 2^{k}-1$ \textit{if and only if}$\,$ $s\left( A\right) \leq
2^{k}-1.$

\textit{ii) If}\thinspace \thinspace \thinspace \underline{$s$}$\left(
A\right) =n\,$ \textit{and\thinspace } $k\geq 2\,$ \textit{such
that\thinspace } $2^{k-1}\leq n<2^{k},$ \textit{then} $\,s\left( A\right)
\leq 2^{k}-1.$

\textit{iii)} \textit{If} \underline{$s$}$\left( A\right) =1$ \textit{then} $%
\,s\left( A\right) \leq 2.\medskip $

\textbf{Proof.} i) For $k\geq 2,$ then \underline{$s$}$\left( A\right) \leq
2^{k}-1$ if and only if $<1>\perp \left( \left( 2^{k}-1\right) \times
T_{P}\right) $\textit{\ }is isotropic. This is equivalent with $s\left(
A\right) \leq 2^{k}-1$.

ii) If\thinspace \thinspace \thinspace \thinspace $n<2^{k},\,\,k\geq 2,$ it
results $n\leq 2^{k}-1,$ and we apply i).

iii) We have that \underline{$s$}$\left( A\right) =1$ if and only
if\thinspace \thinspace $<1>\perp T_{P}=T_{C}\;$\textit{\ }is isotropic or%
\textit{\ }$2\times T_{P}$\textit{\ }is isotropic. If $\,\,2\times T_{P}$%
\textit{\ }is isotropic, then it is universal and represents $-1.$ Therefore 
$\,s\left( A\right) \leq 2. $ If\thinspace \thinspace $T_{C\text{ }}\,\,$is
isotropic, then $T_{P}$ is isotropic, then is universal and represents $-1.$
We obtain $s\left( A\right) =1.\Box \medskip $

\textbf{Proposition 3.9.} \textit{With the above notations, we have}:

\textit{i) For} $k\geq 2,\,$ \textit{if} \thinspace \underline{$s$}$\left(
A\right) =2^{k}-1$ \textit{\thinspace then } $s\left( A\right) =2^{k}-1.$

\textit{ii)} \textit{For} $k\geq 2,\,$ \textit{if} \thinspace $s\left(
A\right) =2^{k}$ \textit{\thinspace then } \underline{$s$}$\left( A\right)
=2^{k}.$

\textit{iii)} \textit{For} $k\geq 1,\,$ \textit{if} \thinspace $s\left(
A\right) =2^{k}+1$ \textit{\thinspace then } $\,\,$\underline{$s$}$\left(
A\right) =2^{k}+1$ \ \textit{or\thinspace \ } \underline{$s$}$\left(
A\right) =2^{k}.\medskip $

\textbf{Proof.} i) From Proposition 3.8., if \thinspace \thinspace 
\underline{$s$}$\left( A\right) =2^{k}-1$ \thinspace then\thinspace
\thinspace $s\left( A\right) \leq 2^{k}-1.$ Since \underline{$s$}$\left(
A\right) \leq s\left( A\right) ,$ therefore $s\left( A\right) =2^{k}-1.$

ii) If \thinspace \underline{$s$}$\left( A\right) \leq 2^{k}-1$ we have $%
\,\,s\left( A\right) \leq 2^{k}-1,$ false$.$

iii) For $k\geq 1,$ if $\,s\left( A\right) =2^{k}+1,$ since \underline{$s$}$%
\left( A\right) \leq s\left( A\right) ,$ we obtain that \underline{$s$}$%
\left( A\right) \leq 2^{k}+1.$ If \underline{$s$}$\left( A\right) \leq
2^{k}-1,$ then $s\left( A\right) \leq 2^{k}-1,$ false.$\Box \medskip $

Some of the above results are proved by \thinspace \thinspace J. O'Shea in
[16] and in [17] for the quaternions and octonions$.$%
\begin{equation*}
\end{equation*}%
\textbf{Acknowledgements.} I would like to thank Susanne Pumpl\"{u}n for her
many encouragements and suggestions which helped me to improve this paper. I
am deeply grateful for all.%
\begin{equation*}
\end{equation*}

\bigskip
\bigskip
\bigskip
\bigskip
\bigskip
\bigskip
\bigskip

{\large References\medskip }

[1] Brown, R. B., \textit{On generalized Cayley-Dickson algebras}, Pacific
J. of Math.\textbf{\ 20(3)}(1967), 415-422 .

[2] Dieterich, E., \textit{Real quadratic division algebras}, Comm. Algebra. 
\textbf{28(2)}(2000), 941-947.

[3] Elman, R., Lam, T. Y., \textit{Pfister forms and K-theory of fields,} J.
Alg. \textbf{23}(1972), 181-213.

[4] Flaut, C., \textit{\ } \textit{Divison algebras with dimension }$%
2^{t},t\in $\textit{\ }$\mathbb{N}$, Analele Stiintifice ale Universitatii
\textquotedblleft Ovidius\textquotedblright\ Constanta, Seria Matematica, 
\textbf{13(2)}(2006), 31-38.

[5] Hoffman, D. W., \textit{Isotropy of quadratic forms over the function
field of a quadric}, Mathematische Zeitschrift. 220(3), 461-476 \ (1995).

[6] Hoffman, D. W., \textit{Levels and sublevels of quaternion algebras,}
Mathematical Proceedings of the Royal Irish Academy. \textbf{110A(1)}(2010),
95-98.

[7] Hoffman, D. W., \textit{Levels of quaternion algebras}, Archiv der
Mathematik. \textbf{90(5)}(2008), 401-411.

[8] Knebusch, M., \textit{Generic splitting of quadratic forms I}, Proc.
London Math. Soc. \textbf{33}(1976), 65-93.

[9] Koprowski, P., \ \textit{Sums of squares of pure quaternions,} Proc.
Roy. Irish Acad. \textbf{98(1)}(1998), 63-65.

[10] K\"{u}skemper, M., \ Wadsworth, A., \textit{A quaternion algebra of
sublevel \ 3}, Bull. Soc. Math. Belg. S\'{e}r. B.\textbf{\ 43(2)}(1991),
181-185.

[11] Laghribi, A., Mammone, P., \textit{On the level of a quaternion algebra}%
, Comm. Algebra. \textbf{29(4)}(2001), 1821-1828.

[12] \thinspace Leep, D. B., \textit{Levels of division algebras}, Glasgow
Math. J. \textbf{32}(1990), 365-370.

[13] Lewis, D. W., \textit{Levels and sublevels of division algebras,} Proc.
Roy. Irish Acad. Sect. A. \textbf{87(1)}(1987), 103-106.

[14] Lewis, D. W., \textit{Levels of quaternion algebras}, Rocky Mountain J.
Math. \textbf{19}(1989), 787-792 .

[15] Osborn, J. M., \textit{Quadratic Division Algebras}, Trans. Amer. Math.
Soc. \textbf{115}(1962), 202-221.

[16] O' Shea, J., \textit{Bounds on the levels of composition algebras},
Mathematical Proceedings of the Royal Irish Academy. \textbf{110A(1)}(2010),
21-30.

[17] \ O' Shea, J., \textit{Levels and sublevels of composition algebras,}
Indag. Mathem. \textbf{18(1)}(2007), 147-159.

[18] Pfister, A., \textit{Zur Darstellung von-I als Summe von quadraten in
einem K\"{o}rper,} J. London Math. Soc. \textbf{40}(1965), 159-165.

[19] Pumpl\"{u}n, S., \textit{Sums of squares in octonion algebras}, Proc.
Amer. Math. Soc. \textbf{133}(2005), 3143-3152.

[20] Schafer, R. D., \textit{An Introduction to Nonassociative Algebras,}
New York. Academic Press, 1966.

[21] Schafer, R. D., \textit{On the algebras formed by the Cayley-Dickson
process,} Amer. J. Math. \textbf{76}(1954), 435-446.

[22] Scharlau, W., \textit{Quadratic and Hermitian Forms},
Berlin-Heidelberg-New York-Tokyo, Springer-Verlag, \ 1985.

\end{document}